\documentclass[english]{article}
\usepackage[T1]{fontenc}
\usepackage[latin9]{inputenc}
\usepackage{geometry}
\usepackage{babel}
\usepackage{url}
\usepackage{bm}
\usepackage{amsmath}
\usepackage{amsthm}
\usepackage{amssymb}
\usepackage[authoryear,sort,comma,round]{natbib}
 
 \usepackage[colorlinks,linkcolor=magenta,citecolor=blue,pagebackref=true]{hyperref}
\renewcommand*{\backrefalt}[4]{%
	\ifcase #1 \footnotesize{(Not cited.)}%
	\or        \footnotesize{(Cited on page~#2.)}%
	\else      \footnotesize{(Cited on pages~#2.)}%
	\fi}

\makeatletter
\newcommand{\lyxaddress}[1]{
	\par {\raggedright #1
	\vspace{1.4em}
	\noindent\par}
}
\theoremstyle{plain}
\newtheorem{thm}{\protect\theoremname}
\theoremstyle{plain}
\newtheorem{lem}{\protect\lemmaname}
\theoremstyle{plain}
\newtheorem{cor}{\protect\corollaryname}

\usepackage{babel}

\makeatother

\providecommand{\corollaryname}{Corollary}
\providecommand{\lemmaname}{Lemma}
\providecommand{\theoremname}{Theorem}

\begin{document}
\title{Approximations of conditional probability density functions in Lebesgue
spaces via mixture of experts models}
\author{Hien Duy Nguyen$^{*1}$ \and TrungTin Nguyen$^{2}$ \and Faicel
Chamroukhi$^{2}$ \and Geoffrey John McLachlan$^{3}$ }
\maketitle

\lyxaddress{$^{1}$Department of Mathematics and Statistics, La Trobe University,
Bundoora Victoria, Australia. \\
$^{2}$Normandie Univ, UNICAEN, CNRS,
LMNO, 14000 Caen, France. \\
$^{3}$School of Mathematics and Physics,
University of Queensland, St. Lucia Brisbane, Australia.\\ $^{*}$Corresponding
author---email: \url{h.nguyen5@latrobe.edu.au}.}
\begin{abstract}
Mixture of experts (MoE) models are widely applied for conditional
probability density estimation problems. We demonstrate the richness
of the class of MoE models by proving denseness results in Lebesgue
spaces, when inputs and outputs variables are both compactly supported.
We further prove an almost uniform convergence result when the input
is univariate. Auxiliary lemmas are proved regarding the richness
of the soft-max gating function class, and their relationships to
the class of Gaussian gating functions. 
\end{abstract}

\section{Introduction}

Mixture of experts (MoE) models are a widely applicable class of conditional
probability density approximations that have been considered as solution
methods across the spectrum of statistical and machine learning problems;
see, for example, the reviews of \citet{Yuksel2012}, \citet{Masoudnia2014},
and \citet{Nguyen2018}.

Let $\mathbb{Z}=\mathbb{X}\times\mathbb{Y}$, where $\mathbb{X}\subseteq\mathbb{R}^{d}$
and $\mathbb{Y}\subseteq\mathbb{R}^{q}$, for $d,q\in\mathbb{N}$.
Suppose that the input and output random variables, $\bm{X}\in\mathbb{X}$
and $\bm{Y}\in\mathbb{Y}$, are related via the conditional probability
density function (PDF) $f\left(\bm{y}|\bm{x}\right)$ in the functional
class:

\[
\mathcal{F}=\left\{ f:\mathbb{Z}\rightarrow\left[0,\infty\right)|\int_{\mathbb{Y}}f\left(\bm{y}|\bm{x}\right)\text{d}\lambda\left(\bm{y}\right)=1,\forall\bm{x}\in\mathbb{X}\right\} \text{,}
\]
where $\lambda$ denotes the Lebesgue measure. The MoE approach seeks
to approximate the unknown target conditional PDF $f$ by a function
of the MoE form:

\[
m\left(\bm{y}|\bm{x}\right)=\sum_{k=1}^{K}\text{Gate}_{k}\left(\bm{x}\right)\text{Expert}_{k}\left(\bm{y}\right)\text{,}
\]
where $\mathbf{Gate}=\left(\text{Gate}_{k}\right)_{k\in\left[K\right]}\in\mathcal{G}^{K}$
($\left[K\right]=\left\{ 1,\dots,K\right\} $), $\text{Expert}_{1},\dots,\text{Expert}_{K}\in\mathcal{E}$,
and $K\in\mathbb{N}$. Here, we say that $m$ is a $K\text{-component}$
MoE model with gates arising from the class $\mathcal{G}^{K}$ and
experts arising from $\mathcal{E}$, where $\mathcal{E}$ is a class
of PDFs with support $\mathbb{Y}$.

The most popular choices for $\mathcal{G}^{K}$ are the parametric
soft-max and Gaussian gating classes:
\[
\mathcal{G}_{S}^{K}=\left\{ \mathbf{Gate}=\left(\text{Gate}_{k}\left(\cdot;\bm{\gamma}\right)\right)_{k\in\left[K\right]}|\forall k\in\left[K\right],\text{Gate}_{k}\left(\cdot;\bm{\gamma}\right)=\frac{\exp\left(a_{k}+\bm{b}_{k}^{\top}\cdot\right)}{\sum_{l=1}^{K}\exp\left(a_{l}+\bm{b}_{l}^{\top}\cdot\right)},\bm{\gamma}\in\mathbb{G}_{S}^{K}\right\} 
\]
and
\[
\mathcal{G}_{G}^{K}=\left\{ \mathbf{Gate}=\left(\text{Gate}_{k}\left(\cdot;\bm{\gamma}\right)\right)_{k\in\left[K\right]}|\forall k\in\left[K\right],\text{Gate}_{k}\left(\cdot;\bm{\gamma}\right)=\frac{\pi_{k}\phi\left(\cdot;\bm{\nu}_{k},\bm{\Sigma}_{k}\right)}{\sum_{l=1}^{K}\pi_{l}\phi\left(\cdot;\bm{\nu}_{l},\bm{\Sigma}_{l}\right)},\bm{\gamma}\in\mathbb{G}_{G}^{K}\right\} \text{,}
\]
respectively, where
\[
\mathbb{G}_{S}^{K}=\left\{ \bm{\gamma}=\left(a_{1},\dots,a_{K},\bm{b}_{1},\dots,\bm{b}_{K}\right)\in\mathbb{R}^{K}\times\left(\mathbb{R}^{d}\right)^{K}\right\} 
\]
and
\[
\mathbb{G}_{G}^{K}=\left\{ \bm{\gamma}=\left(\bm{\pi},\bm{\nu}_{1},\dots,\bm{\nu}_{K},\bm{\Sigma}_{1},\dots,\bm{\Sigma}_{K}\right)\in\Pi_{K-1}\times\left(\mathbb{R}^{d}\right)^{K}\times\mathbb{S}_{d}^{K}\right\} \text{.}
\]
Here,
\[
\phi\left(\cdot;\bm{\nu},\bm{\Sigma}\right)=\left|2\pi\bm{\Sigma}\right|^{-1/2}\exp\left[-\frac{1}{2}\left(\cdot-\bm{\nu}\right)^{\top}\bm{\Sigma}^{-1}\left(\cdot-\bm{\nu}\right)\right]
\]
is the multivariate normal density function with mean vector $\bm{\nu}$
and covariance matrix $\bm{\Sigma}$, $\bm{\pi}^{\top}=\left(\pi_{1},\dots,\pi_{K}\right)$
is a vector of weights in the simplex:
\[
\Pi_{K-1}=\left\{ \bm{\pi}=\left(\pi_{k}\right)_{k\in\left[K\right]}|\forall k\in\left[K\right],\text{\ensuremath{\pi_{k}>0}},\sum_{k=1}^{K}\pi_{k}=1\right\} \text{,}
\]
and $\mathbb{S}_{d}$ is the class of $d\times d$ symmetric positive
definite matrices. The soft-max and Gaussian gating classes were first
introduced by \citet{Jacobs1991} and \citet{Jordan1995}, respectively.
Typically, one chooses experts that arise from some location-scale
class:

\[
\mathcal{E}_{\psi}=\left\{ g_{\psi}\left(\cdot;\bm{\mu},\sigma\right):\mathbb{Y}\rightarrow\left[0,\infty\right)|g_{\psi}\left(\cdot;\bm{\mu},\sigma\right)=\frac{1}{\sigma^{q}}\psi\left(\frac{\cdot-\bm{\mu}}{\sigma}\right),\bm{\mu}\in\mathbb{R}^{q},\sigma\in\left(0,\infty\right)\right\} \text{,}
\]
where $\psi$ is a PDF, with respect to $\mathbb{R}^{q}$ in the sense
that $\psi:\mathbb{R}^{q}\rightarrow\left[0,\infty\right)$ and $\int_{\mathbb{R}^{q}}\psi\left(\bm{y}\right)\text{d}\lambda\left(\bm{y}\right)=1$.

We shall say that $f\in\mathcal{L}_{p}\left(\mathbb{Z}\right)$ for
any $p\in\left[1,\infty\right)$ if
\[
\left\Vert f\right\Vert _{p,\mathbb{Z}}=\left(\int_{\mathbb{Z}}\left|\mathbf{1}_{\mathbb{Z}}f\right|^{p}\text{d}\lambda\left(\bm{z}\right)\right)^{1/p}<\infty\text{,}
\]
where $\mathbf{1}_{\mathbb{Z}}$ is the indicator function that takes
value $1$ when $\bm{z}\in\mathbb{Z}$, and 0 otherwise. Further,
we say that $f\in\mathcal{L}_{\infty}\left(\mathbb{Z}\right)$ if
\[
\left\Vert f\right\Vert _{\infty,\mathbb{Z}}=\inf\left\{ a\ge0|\lambda\left(\left\{ \bm{z}\in\mathbb{Z}|\left|f\left(\bm{z}\right)\right|>a\right\} \right)=0\right\} <\infty\text{.}
\]
We shall refer to $\left\Vert \cdot\right\Vert _{p,\mathbb{Z}}$ as
the $\mathcal{L}_{p}$ norm on $\mathbb{Z}$, for $p\in\left[0,\infty\right]$,
and where the context is obvious, we shall drop the reference to $\mathbb{Z}$.

Suppose that the target conditional PDF $f$ is in the class $\mathcal{F}_{p}=\mathcal{F}\cap\mathcal{L}_{p}$.
We address the problem of approximating $f$, with respect to the
$\mathcal{L}_{p}$ norm, using MoE models in the soft-max and Gaussian
gated classes, 
\begin{eqnarray*}
\mathcal{M}_{S}^{\psi} & = & \biggl\{ m_{K}^{\psi}:\mathbb{Z}\rightarrow\left[0,\infty\right)|m_{K}^{\psi}\left(\bm{y}|\bm{x}\right)=\sum_{k=1}^{K}\text{Gate}_{k}\left(\bm{x}\right)g_{\psi}\left(\bm{y};\bm{\mu}_{k},\sigma_{k}\right)\text{,}\\
 &  & \quad g_{\psi}\in\mathcal{E}_{\psi}\cap\mathcal{L}_{\infty},\mathbf{Gate}\in\mathcal{G}_{S}^{K},\bm{\mu}_{k}\in\mathbb{Y},\sigma_{k}\in\left(0,\infty\right),k\in\left[K\right],K\in\mathbb{N}\biggr\}
\end{eqnarray*}
and 
\begin{eqnarray*}
\mathcal{M}_{G}^{\psi} & = & \biggl\{ m_{K}^{\psi}:\mathbb{Z}\rightarrow\left[0,\infty\right)|m_{K}^{\psi}\left(\bm{y}|\bm{x}\right)=\sum_{k=1}^{K}\text{Gate}_{k}\left(\bm{x}\right)g_{\psi}\left(\bm{y};\bm{\mu}_{k},\sigma_{k}\right)\text{,}\\
 &  & \quad g_{\psi}\in\mathcal{E}_{\psi}\cap\mathcal{L}_{\infty},\mathbf{Gate}\in\mathcal{G}_{G}^{K},\bm{\mu}_{k}\in\mathbb{Y},\sigma_{k}\in\left(0,\infty\right),k\in\left[K\right],K\in\mathbb{N}\biggr\}\text{,}
\end{eqnarray*}
by showing that both $\mathcal{M}_{S}^{\psi}$ and $\mathcal{M}_{G}^{\psi}$
are dense in the class $\mathcal{F}_{p}$, when $\mathbb{X}=\left[0,1\right]^{d}$
and $\mathbb{Y}$ is a compact subset of $\mathbb{R}^{q}$. Our denseness
results are enabled by the indicator function approximation result
of \citet{Jiang1999}, and the finite mixture model denseness theorems
of \citet{Nguyen:2020aa} and \citet{Nguyen:2020ab}.

Our theorems contribute to an enduring continuity of sustained interest
in the approximation capabilities of MoE models. Related to our results
are contributions regarding the approximation capabilities of the
conditional expectation function of the classes $\mathcal{M}_{S}^{\psi}$
and $\mathcal{M}_{G}^{\psi}$ \citep{Wang1992,Zeevi1998,Jiang1999,Krzyzak:2005aa,MendesJiang2012,Nguyen2016a,Nguyen:2019aa}
and the approximation capabilities of subclasses of $\mathcal{M}_{S}^{\psi}$
and $\mathcal{M}_{G}^{\psi}$, with respect to the Kullback--Leibler
divergence \citep{Jiang1999a,Norets2010,Norets2014}. Our results
can be seen as complements to the Kullback--Leibler approximation
theorems of \citet{Norets2010} and \citet{Norets2014}, by the relationship
between the Kullback--Leibler divergence and the $\mathcal{L}_{2}$
norm \citep{ZeeviMeir1997}. That is, when $f>1/\kappa$, for all
$\left(\bm{x},\bm{y}\right)\in\mathbb{Z}$ and some constant $\kappa>0$,
we have that the integrated conditional Kullback--Leibler
divergence considered by \citet{Norets2014}:
\[
\int_{\mathbb{X}}D\left(f\left(\cdot|\bm{x}\right)\Vert m_{K}^{\psi}\left(\cdot|\bm{x}\right)\right)\text{d}\lambda\left(\bm{x}\right)=\int_{\mathbb{X}}\int_{\mathbb{Y}}f\left(\bm{y}|\bm{x}\right)\log\frac{f\left(\bm{y}|\bm{x}\right)}{m_{K}^{\psi}\left(\bm{y}|\bm{x}\right)}\text{d}\lambda\left(\bm{y}\right)\text{d}\lambda\left(\bm{x}\right)
\]
satisfies
\[
\int_{\mathbb{X}}D\left(f\left(\cdot|\bm{x}\right)\Vert m_{K}^{\psi}\left(\cdot|\bm{x}\right)\right)\text{d}\lambda\left(\bm{x}\right)\le\kappa^{2}\left\Vert f-m_{K}^{\psi}\right\Vert _{2,\mathbb{Z}}^{2}\text{,}
\]
and thus a good approximation in the integrated Kullback--Leibler
divergence is guaranteed if one can find a good approximation in the
$\mathcal{L}_{2}$ norm, which is guaranteed by our main result.

The remainder of the manuscript proceeds as follows. The main result
is presented in Section \ref{sec:Main-results}. Technical lemmas
are provided in Section \ref{sec:Technical-lemmas}. The proofs of
our results are then presented in Section \ref{sec:Proofs-of-main}.
Proofs of required lemmas that do not appear elsewhere are provided
in Section \ref{sec:Proofs-of-lemmas}. A summary of our work and
some conclusions are drawn in Section \ref{sec:Summary-and-conclusions}.

\section{\label{sec:Main-results}Main results}

Denote the class of bounded functions on $\mathbb{Z}$ by 
\[
\mathcal{B}\left(\mathbb{Z}\right)=\left\{ f\in\mathcal{L}_{\infty}\left(\mathbb{Z}\right)|\exists a\in\left[0,\infty\right)\text{, such that}\text{ }\left|f\left(\bm{z}\right)\right|\le a\text{, }\forall\bm{z}\in\mathbb{Z}\right\} \text{,}
\]
and write its norm as $\left\Vert f\right\Vert _{\mathcal{B}\left(\mathbb{Z}\right)}=\sup_{\bm{z}\in\mathbb{Z}}\left|f\left(\bm{z}\right)\right|$.
Further, let $\mathcal{C}$ denote the class of continuous functions.
Note that if $\mathbb{Z}$ is compact and $f\in\mathcal{C}$,
then $f\in\mathcal{B}$. 
\begin{thm}
\label{thm main} Assume that $\mathbb{X}=\left[0,1\right]^{d}$
for $d\in\mathbb{N}$. There exists a sequence $\left\{ m_{K}^{\psi}\right\} _{K\in\mathbb{N}}\subset\mathcal{M}_{S}^{\psi}$,
such that if $\mathbb{Y}\subset\mathbb{R}^{q}$ is compact, $f\in\mathcal{F}\cap\mathcal{C}$,
and $\psi\in\mathcal{C}\left(\mathbb{R}^{q}\right)$ is a PDF on support
$\mathbb{R}^{q}$, then $\lim_{K\rightarrow\infty}\left\Vert f-m_{K}^{\psi}\right\Vert _{p}=0$,
for $p\in\left[1,\infty\right)$ .
\end{thm}
Since convergence in Lebesgue spaces does not imply point-wise modes
of convergence, the following result is also useful and interesting
in some restricted scenarios. Here, we note that the mode of convergence
is almost uniform, which implies almost everywhere convergence and
convergence in measure (cf. \citealt[Lem 7.10 and Thm. 7.11]{Bartle:1995aa}).
The almost uniform convergence of $\left\{ m_{K}^{\psi}\right\} _{K\in\mathbb{N}}$
to $f$ in the following result is to be understood in the sense of
\citet[Def. 7.9]{Bartle:1995aa}. That is, for every $\delta>0$,
there exists a set $\mathbb{E}_{\delta}\subset\mathbb{Z}$ with $\lambda\left(\mathbb{Z}\right)<\delta$,
such that $\left\{ m_{K}^{\psi}\right\} _{K\in\mathbb{N}}$ converges
to $f$, uniformly on $\mathbb{Z}\backslash\mathbb{E}_{\delta}$. 
\begin{thm}
\label{thm main d=00003D00003D1} Assume that $\mathbb{X}=\left[0,1\right]$.
There exists a sequence $\left\{ m_{K}^{\psi}\right\} _{K\in\mathbb{N}}\subset\mathcal{M}_{S}^{\psi}$,
such that if $\mathbb{Y}\subset\mathbb{R}^{q}$ is compact, $f\in\mathcal{F}\cap\mathcal{C}$,
and $\psi\in\mathcal{C}\left(\mathbb{R}^{q}\right)$ is a PDF on support
$\mathbb{R}^{q}$, then $\lim_{K\rightarrow\infty}m_{K}^{\psi}=f$,
almost uniformly. 
\end{thm}
The following result establishes the connection between the gating
classes $\mathcal{G}_{S}^{K}$ and $\mathcal{G}_{G}^{K}$. 
\begin{lem}
\label{lem equiv}For each $K\in\mathbb{N}$, $\mathcal{G}_{S}^{K}\subset\mathcal{G}_{G}^{K}$.
Further, if we define the class of Gaussian gating vectors with equal
covariance matrices: 
\[
\mathcal{G}_{E}^{K}=\left\{ \mathbf{Gate}=\left(\mathrm{Gate}_{k}\left(\cdot;\bm{\gamma}\right)\right)_{k\in\left[K\right]}|\forall k\in\left[K\right],\mathrm{Gate}_{k}\left(\cdot;\bm{\gamma}\right)=\frac{\pi_{k}\phi\left(\cdot;\bm{\nu}_{k},\bm{\Sigma}\right)}{\sum_{l=1}^{K}\pi_{l}\phi\left(\cdot;\bm{\nu}_{l},\bm{\Sigma}\right)},\bm{\gamma}\in\mathbb{G}_{E}^{K}\right\} \text{,}
\]
where 
\[
\mathbb{G}_{E}^{K}=\left\{ \bm{\gamma}=\left(\bm{\pi},\bm{\nu}_{1},\dots,\bm{\nu}_{K},\bm{\Sigma}\right)\in\Pi_{K-1}\times\left(\mathbb{R}^{d}\right)^{K}\times\mathbb{S}_{d}\right\} \text{,}
\]
then $\mathcal{G}_{E}^{K}\subset\mathcal{G}_{S}^{K}$. 
\end{lem}
We can directly apply Lemma \ref{lem equiv} to establish the following
corollary to Theorems \ref{thm main} and \ref{thm main d=00003D00003D1},
regarding the approximation capability of the class $\mathcal{M}_{G}^{\psi}$. 
\begin{cor}
\label{cor: Gaussian gate}Theorems \ref{thm main} and \ref{thm main d=00003D00003D1}
hold when $\mathcal{M}_{S}^{\psi}$ is replaced by $\mathcal{M}_{G}^{\psi}$
in their statements. 
\end{cor}

\section{\label{sec:Technical-lemmas}Technical lemmas}

Let $\mathbb{K}^{n}=\left\{ \left(k_{1},\dots,k_{d}\right)\in\left[n\right]^{d}\right\} $
and $\kappa:\mathbb{K}^{n}\rightarrow\left[n^{d}\right]$ be a bijection 
for each $n\in\mathbb{N}$. For each $\left(k_{1},\dots,k_{d}\right)\in\mathbb{K}^{n}$
and $k\in\left[n^{d}\right]$, we define $\mathbb{X}_{k}^{n}=\mathbb{X}_{\kappa\left(k_{1},\dots,k_{d}\right)}^{n}=\prod_{i=1}^{d}\mathbb{I}_{k_{i}}^{n}$,
where $\mathbb{I}_{k_{i}}^{n}=\left[\left(k_{i}-1\right)/n,k_{i}/n\right)$
for $k_{i}\in\left[n-1\right]$, and $\mathbb{I}_{n}^{n}=\left[\left(n-1\right)/n,1\right]$.

We call $\left\{ \mathbb{X}_{k}^{n}\right\} _{k\in\left[n^{d}\right]}$
a fine partition of $\mathbb{X}$, in the sense that $\mathbb{X}=\left[0,1\right]^{d}=\bigcup_{k=1}^{n^{d}}\mathbb{X}_{k}^{n}$,
for each $n$, and that $\lambda\left(\mathbb{X}_{k}^{n}\right)=n^{-d}$
gets smaller, as $n$ increases. The following result from \citet{Jiang1999}
establishes the approximation capability of soft-max gates. 
\begin{lem}[Jiang and Tanner, 1999, p. 1189]
\label{Lem Jiang}For each $n\in\mathbb{N}$, $p\in\left[1,\infty\right)$
and $\epsilon>0$, there exists a gating functions 
\[
\mathbf{Gate}=\left(\mathrm{Gate}_{k}\left(\cdot;\bm{\gamma}\right)\right)_{k\in\left[n^{d}\right]}\in\mathcal{G}_{S}^{n^{d}}
\]
for some $\bm{\gamma}\in\mathbb{G}_{S}^{n^{d}}$, such that 
\[
\sup_{k\in\left[n^{d}\right]}\left\Vert \mathbf{1}_{\left\{ \bm{x}\in\mathbb{X}_{k}^{n}\right\} }-\mathrm{Gate}_{k}\left(\cdot;\bm{\gamma}\right)\right\Vert _{p,\mathbb{X}}\le\epsilon\text{.}
\]
\end{lem}
When, $d=1$, we have also the following almost uniform convergence
alternative to Lemma \ref{Lem Jiang}. 
\begin{lem}
\label{lem almost uniform jiang}Let $\mathbb{X}=\left[0,1\right]$.
Then, for each $n\in\mathbb{N}$, there exists a sequence of gating
functions:
\[
\left\{ \mathbf{Gate}_{l}=\left(\mathrm{Gate}_{k}\left(\cdot;\bm{\gamma}_{l}\right)\right)_{k\in\left[n^{d}\right]}\right\} _{l\in\mathbb{N}}\subset\mathcal{G}_{S}^{n}\text{,}
\]
defined by $\left\{ \bm{\gamma}_{l}\right\} _{l\in\mathbb{N}}\subset\mathbb{G}_{S}^{n}$,
such that 
\[
\mathrm{Gate}_{k}\left(\cdot;\bm{\gamma}_{l}\right)\rightarrow\mathbf{1}_{\left\{ \bm{x}\in\mathbb{X}_{k}^{n}\right\} }\text{,}
\]
almost uniformly, simultaneously for all $k\in\left[n^{d}\right]$. 
\end{lem}
For PDF $\psi$ on support $\mathbb{R}^{q}$, define the class of
finite mixture models by
\begin{eqnarray*}
\mathcal{H}^{\psi} & = & \biggl\{ h_{K}^{\psi}:\mathbb{R}^{q}\rightarrow\left[0,\infty\right)|h_{K}^{\psi}\left(\bm{y}\right)=\sum_{k=1}^{K}c_{k}g_{\psi}\left(\bm{y};\bm{\mu}_{k},\sigma_{k}\right)\text{,}\\
 &  & \quad g_{\psi}\in\mathcal{E}_{\psi}\cap\mathcal{L}_{\infty},\left(c_{k}\right)_{k\in\left[K\right]}\in\Pi_{K-1},\bm{\mu}_{k}\in\mathbb{Y},\sigma_{k}\in\left(0,\infty\right),k\in\left[K\right],K\in\mathbb{N}\biggr\}\text{.}
\end{eqnarray*}
We require the following result, from \citet{Nguyen:2020aa}, regarding the approximation capabilities
of $\mathcal{H}^{\psi}$. 
\begin{lem}[Nguyen et al., 2020a, Thm. 2(b)]
\label{lem tin lem} If $f\in\mathcal{C}\left(\mathbb{Y}\right)$
is a PDF on $\mathbb{Y}$, $\psi\in\mathcal{C}\left(\mathbb{R}^{q}\right)$
is a PDF on $\mathbb{R}^{q}$, and $\mathbb{Y}\subset\mathbb{R}^{q}$
is compact, then there exists a sequence $\left\{ h_{K}^{\psi}\right\} _{K\in\mathbb{N}}\subset\mathcal{H}^{\psi}$,
such that $\lim_{K\rightarrow\infty}\left\Vert f-h_{K}^{\psi}\right\Vert _{\mathcal{B}\left(\mathbb{Y}\right)}=0$. 
\end{lem}

\section{\label{sec:Proofs-of-main}Proofs of main results}

\subsection{Proof of Theorem \ref{thm main}}

To prove the result, it suffices to show that for each $\epsilon>0$,
there exists a $m_{K}^{\psi}\in\mathcal{M}_{S}^{\psi}$, such that
\[
\left\Vert f-m_{K}^{\psi}\right\Vert _{p}<\epsilon\text{.}
\]

The main steps of the proof are as follows. We firstly approximate
$f\left(\bm{y}|\bm{x}\right)$ by 
\begin{equation}
\upsilon_{n}\left(\bm{y}|\bm{x}\right)=\sum_{k=1}^{n^{d}}\mathbf{1}_{\left\{ \bm{x}\in\mathbb{X}_{k}^{n}\right\} }f\left(\bm{y}|\bm{x}_{k}^{n}\right)\text{,}\label{eq: discrete fine approx}
\end{equation}
where $\bm{x}_{k}^{n}\in\mathbb{X}_{k}^{n}$, for each $k\in\left[n^{d}\right]$,
such that 
\begin{equation}
\left\Vert f-\upsilon_{n}\right\Vert _{p}<\frac{\epsilon}{3}\text{,}\label{eq: res 1}
\end{equation}
for all $n\ge N_{1}\left(\epsilon\right)$, for some sufficiently
large $N_{1}\left(\epsilon\right)\in\mathbb{N}$. Then we approximate
$\upsilon_{n}\left(\bm{y}|\bm{x}\right)$ by

\begin{equation}
\eta_{n}\left(\bm{y}|\bm{x}\right)=\sum_{k=1}^{n^{d}}\text{Gate}_{k}\left(\bm{x};\bm{\gamma}_{n}\right)f\left(\bm{y}|\bm{x}_{k}^{n}\right)\text{,}\label{eq: eta fun}
\end{equation}
where $\bm{\gamma}_{n}\in\mathbb{G}_{S}^{n^{d}}$ and $\mathbf{Gate}=\left(\mathrm{Gate}_{k}\left(\cdot;\bm{\gamma}_{n}\right)\right)_{k\in\left[n^{d}\right]}\in\mathcal{G}_{S}^{n^{d}}$,
so that 
\begin{align}
\left\Vert \upsilon_{n}-\eta_{n}\right\Vert _{p} & \le\sup_{k\in\left[n^{d}\right]}\left\Vert \mathrm{Gate}_{k}\left(\cdot;\bm{\gamma}\right)-\mathbf{1}_{\left\{ \bm{x}\in\mathbb{X}_{k}^{n}\right\} }\right\Vert _{p,\mathbb{X}}\sum_{k=1}^{n^{d}}\left\Vert f\left(\cdot|\bm{x}_{k}^{n}\right)\right\Vert _{p,\mathbb{Y}}<\frac{\epsilon}{3}\text{,}\label{eq: res 2}
\end{align}
using Lemma \ref{Lem Jiang}.

Finally, we approximate $\eta_{n}\left(\bm{y}|\bm{x}\right)$ by $m_{K_{n}}^{\psi}\left(\bm{y}|\bm{x}\right)$,
where 
\begin{equation}
m_{K_{n}}^{\psi}\left(\bm{y}|\bm{x}\right)=\sum_{k=1}^{n^{d}}\text{Gate}_{k}\left(\bm{x};\bm{\gamma}\right)h_{n_{k}}^{k}\left(\bm{y}|\bm{x}_{k}^{n}\right)\label{eq: all together approx}
\end{equation}
and 
\begin{equation}
h_{n_{k}}^{k}\left(\bm{y}|\bm{x}_{k}^{n}\right)=\sum_{i=1}^{n_{k}}c_{i}^{k}g_{\psi}\left(\bm{y};\bm{\mu}_{i}^{k},\sigma_{i}^{k}\right)\in\mathcal{H}^{\psi}\label{eq: individual approx}
\end{equation}
for $n_{k}\in\mathbb{N}$ ($k\in\left[n^{d}\right]$), such that $K_{n}=\sum_{k=1}^{n^{d}}n_{k}$.
Here, we establish that there exists $N_{2}\left(\epsilon,n,\bm{\gamma}_{n}\right)\in\mathbb{N}$,
so that when $n_{k}\ge N_{2}\left(\epsilon,n,\bm{\gamma}_{n}\right)$,
\begin{equation}
\left\Vert \eta_{n}-m_{K_{n}}^{\psi}\right\Vert _{p}\le\sup_{k\in\left[n^{d}\right]}\left\Vert \mathrm{Gate}_{k}\left(\cdot;\bm{\gamma}\right)\right\Vert _{p,\mathbb{X}}\sum_{k=1}^{n^{d}}\left\Vert f\left(\cdot|\bm{x}_{k}^{n}\right)-h_{n_{k}}^{k}\left(\cdot|\bm{x}_{k}^{n}\right)\right\Vert _{p,\mathbb{Y}}<\frac{\epsilon}{3}\text{.}\label{eq: res 3}
\end{equation}

Results (\ref{eq: res 1})--(\ref{eq: res 3}) then imply that for
each $\epsilon>0$, there exists $N_{1}\left(\epsilon\right)$, $\bm{\gamma}_{n}$,
and $N_{2}\left(\mathbb{\epsilon},n,\bm{\gamma}_{n}\right)$, such
that for all $K_{n}=\sum_{k=1}^{n^{d}}n_{k}$, where $n_{k}\ge N_{2}\left(\epsilon,n,\bm{\gamma}_{n}\right)$
(for each $k\in\left[n^{d}\right]$) and $n\ge N_{1}\left(\epsilon\right)$.
The following inequality results from an application of the triangle inequality: 
\begin{align*}
\left\Vert f-m_{K_{n}}^{\psi}\right\Vert _{p} & \le\left\Vert f-\upsilon_{n}\right\Vert _{p}+\left\Vert \upsilon_{n}-\eta_{n}\right\Vert _{p}+\left\Vert \eta_{n}-m_{K_{n}}^{\psi}\right\Vert _{p}<3\times\frac{\epsilon}{3}=\epsilon\text{.}
\end{align*}

We now focus our attention to proving each of the results: (\ref{eq: res 1})--(\ref{eq: res 3}).
To prove (\ref{eq: res 1}), we note that since $f$ is uniformly
continuous (because $\mathbb{Z}=\mathbb{X}\times\mathbb{Y}$ is compact,
and $f\in\mathcal{C}$), there exists a function (\ref{eq: discrete fine approx}) 
such that for all $\varepsilon>0$,

\begin{equation}
\sup_{\left(\bm{x},\bm{y}\right)\in\mathbb{Z}}\left|f\left(\bm{y}|\bm{x}\right)-\upsilon\left(\bm{y}|\bm{x}\right)\right|<\varepsilon\text{.}\label{eq: pointwise}
\end{equation}

We can construct such an approximation by considering the fact that
as $n$ increases, the diameter $\delta_{n}=\sup_{k\in n^{d}}\mathrm{diam}\left(\mathbb{X}_{k}^{n}\right)$
of the fine partition goes to zero. By the uniform continuity of $f$,
for every $\varepsilon>0$, there exists a $\delta\left(\epsilon\right)>0$,
such that if $\left\Vert \left(\bm{x}_{1},\bm{y}_{1}\right)-\left(\bm{x}_{2},\bm{y}_{2}\right)\right\Vert <\delta\left(\epsilon\right)$,
then $\left|f\left(\bm{y}_{1}|\bm{x}_{1}\right)-f\left(\bm{y}_{2}|\bm{x}_{2}\right)\right|<\varepsilon$,
for all pairs $\left(\bm{x}_{1},\bm{y}_{1}\right),\left(\bm{x}_{2},\bm{y}_{2}\right)\in\mathbb{Z}$.
Here, $\left\Vert \cdot\right\Vert $ denotes the Euclidean norm.
Furthermore, for any $\left(\bm{x},\bm{y}\right)\in\mathbb{Z}$, we
have

\begin{align}
\left|f\left(\bm{y}|\bm{x}\right)-\upsilon_{n}\left(\bm{y}|\bm{x}\right)\right| & =\left|\sum_{k=1}^{n^{d}}\mathbf{1}_{\left\{ \bm{x}\in\mathbb{X}_{k}^{n}\right\} }\left[f\left(\bm{y}|\bm{x}\right)-f\left(\bm{y}|\bm{x}_{k}^{n}\right)\right]\right|\nonumber \\
 & \le\sum_{k=1}^{n^{d}}\mathbf{1}_{\left\{ \bm{x}\in\mathbb{X}_{k}^{n}\right\} }\left|f\left(\bm{y}|\bm{x}\right)-f\left(\bm{y}|\bm{x}_{k}^{n}\right)\right|\text{,}\label{eq: partition}
\end{align}
by the triangle inequality.

Since $\bm{x}_{k}^{n}\in\mathbb{X}_{k}^{n}$, for each $k$ and $n$,
we have the fact that $\left\Vert \left(\bm{x},\bm{y}\right)-\left(\bm{x}_{k}^{n},\bm{y}\right)\right\Vert <\delta_{n}$
for $\left(\bm{x},\bm{y}\right)\in\mathbb{X}_{k}^{n}\times\mathbb{Y}$.
By uniform continuity, for each $\varepsilon$, we can find a sufficiently
small $\delta\left(\epsilon\right)$, such that $\left|f\left(\bm{y}|\bm{x}\right)-f\left(\bm{y}|\bm{x}_{k}^{n}\right)\right|<\varepsilon$,
if $\left\Vert \left(\bm{x},\bm{y}\right)-\left(\bm{x}_{k}^{n},\bm{y}\right)\right\Vert <\delta\left(\epsilon\right)$,
for all $k$. The desired result (\ref{eq: pointwise}) can be obtained
by noting that the right hand side of (\ref{eq: partition}) consists
of only one non-zero summand for any $\left(\bm{x},\bm{y}\right)\in\mathbb{Z}$,
and by choosing $n\in\mathbb{N}$ sufficiently large, so that $\delta_{n}<\delta\left(\epsilon\right)$.

By (\ref{eq: pointwise}), we have the fact that $\upsilon_{n}\rightarrow f$,
point-wise. We can bound $\upsilon_{n}$ as follows:
\begin{equation}
\upsilon_{n}\left(\bm{y}|\bm{x}\right)\le\sum_{i=1}^{n^{p}}\mathbf{1}_{\left\{ \bm{x}\in\mathbb{X}_{k}^{n}\right\} }\sup_{\bm{\zeta}\in\mathbb{Y},\bm{\xi}\in\mathbb{X}}f\left(\bm{\zeta}|\bm{\xi}\right)=\sup_{\bm{\zeta}\in\mathbb{Y},\bm{\xi}\in\mathbb{X}}f\left(\bm{\zeta}|\bm{\xi}\right)\text{,}\label{eq: control v_n}
\end{equation}
where the right-hand side is a constant and is therefore in $\mathcal{L}_{p}$,
since $\mathbb{Z}$ is compact. An application of the Lebesgue dominated
convergence theorem in $\mathcal{L}_{p}$ then yields (\ref{eq: res 1}).

Next we write 
\begin{align*}
\left\Vert \upsilon_{n}-\eta_{n}\right\Vert _{p} & =\left\Vert \sum_{k=1}^{n^{d}}\mathbf{1}_{\left\{ \bm{x}\in\mathbb{X}_{k}^{n}\right\} }f\left(\bm{y}|\bm{x}_{k}^{n}\right)-\sum_{k=1}^{n^{d}}\mathrm{Gate}_{k}\left(\bm{x};\bm{\gamma}_{n}\right)f\left(\bm{y}|\bm{x}_{k}^{n}\right)\right\Vert _{p}\\
 & \le\sum_{k=1}^{n^{d}}\left\Vert \left[\mathbf{1}_{\left\{ \bm{x}\in\mathbb{X}_{k}^{n}\right\} }-\mathrm{Gate}_{k}\left(\bm{x};\bm{\gamma}_{n}\right)\right]f\left(\bm{y}|\bm{x}_{k}^{n}\right)\right\Vert _{p}\text{.}
\end{align*}
Since the norm arguments are separable in $\bm{x}$ and $\bm{y}$,
we apply Fubini's theorem to get 
\begin{align*}
\left\Vert \upsilon_{n}-\eta_{n}\right\Vert _{p} & =\sum_{k=1}^{n^{d}}\left\Vert \left[\mathbf{1}_{\left\{ \bm{x}\in\mathbb{X}_{k}^{n}\right\} }-\mathrm{Gate}_{k}\left(\bm{x};\bm{\gamma}_{n}\right)\right]\right\Vert _{p,\mathbb{X}}\left\Vert f\left(\bm{y}|\bm{x}_{k}^{n}\right)\right\Vert _{p,\mathbb{Y}}\\
 & \le\sup_{k\in\left[n^{d}\right]}\left\Vert \left[\mathbf{1}_{\left\{ \bm{x}\in\mathbb{X}_{k}^{n}\right\} }-\mathrm{Gate}_{k}\left(\bm{x};\bm{\gamma}_{n}\right)\right]\right\Vert _{p,\mathbb{X}}\sum_{k=1}^{n^{d}}\left\Vert f\left(\bm{y}|\bm{x}_{k}^{n}\right)\right\Vert _{p,\mathbb{Y}}
\end{align*}

Because $f\in\mathcal{B}$ and $n^{d}$ is finite, for any fixed $n\in\mathbb{N}$,
we have $C_{1}\left(n\right)=\sum_{k=1}^{n^{d}}\left\Vert f\left(\bm{y}|\bm{x}_{k}^{n}\right)\right\Vert _{p,\mathbb{Y}}<\infty$.
For each $\epsilon>0$, we need to choose a $\bm{\gamma}_{n}\in\mathbb{G}_{S}^{n^{d}}$,
such that

\[
\sup_{k\in\left[n^{d}\right]}\left\Vert \left[\mathbf{1}_{\left\{ \bm{x}\in\mathbb{X}_{k}^{n}\right\} }-\mathrm{Gate}_{k}\left(\bm{x};\bm{\gamma}_{n}\right)\right]\right\Vert _{p,\mathbb{X}}<\frac{\epsilon}{3C_{1}\left(n\right)}\text{,}
\]
which can be achieved via a direct application of Lemma \ref{Lem Jiang}.
We have thus shown (\ref{eq: res 2}).

Lastly, we are required to approximate $f\left(\bm{y}|\bm{x}_{k}^{n}\right)$
for each $k\in\left[n^{d}\right]$, by a function of form (\ref{eq: individual approx}).
Since $\mathbb{Y}$ is compact and $f$ and $\psi$ are continuous,
we can apply of Lemma \ref{lem tin lem}, directly. Note that over
a set of finite measure, convergence in $\left\Vert \cdot\right\Vert _{\mathcal{B}}$
implies convergence in $\mathcal{L}_{p}$ norm, for all $p\in\left[1,\infty\right]$
(cf. \citealt[Prop. 3.9.3]{Oden2010}).

We can then write (\ref{eq: all together approx}) as 
\begin{align}
m_{K_{n}}^{\psi}\left(\bm{y}|\bm{x}\right) & =\sum_{k=1}^{n^{d}}\frac{\exp\left(a_{n,k}+\bm{b}_{n,k}^{\top}\bm{x}\right)}{\sum_{l=1}^{n^{d}}\exp\left(a_{n,l}+\bm{b}_{n,l}^{\top}\bm{x}\right)}h_{n_{k}}^{k}\left(\bm{y}|\bm{x}_{k}^{n}\right)\nonumber \\
 & =\sum_{k=1}^{n^{d}}\sum_{i=1}^{n_{k}}\frac{\exp\left(a_{n,k}+\bm{b}_{n,k}^{\top}\bm{x}\right)}{\sum_{l=1}^{n^{d}}\exp\left(a_{n,l}+\bm{b}_{n,l}^{\top}\bm{x}\right)}\frac{c_{i}^{k}}{\sum_{l=1}^{n_{k}}c_{l}^{k}}g_{\psi}\left(\bm{y};\bm{\mu}_{i}^{k},\sigma_{i}^{k}\right)\nonumber \\
 & =\sum_{k=1}^{n^{d}}\sum_{i=1}^{n_{k}}\frac{\exp\left(\log c_{i}^{k}+a_{n,k}+\bm{b}_{n,k}^{\top}\bm{x}\right)}{\sum_{l=1}^{n^{d}}\sum_{j=1}^{n_{k}}\exp\left(\log c_{j}^{k}+a_{n,l}+\bm{b}_{n,l}^{\top}\bm{x}\right)}g_{\psi}\left(\bm{y};\bm{\mu}_{i}^{k},\sigma_{i}^{k}\right)\text{,}\label{eq: alltogether expand}
\end{align}
where $\bm{\gamma}_{n}=\left(a_{n,1},\dots,a_{n,n^{d}},\bm{b}_{n,1},\dots,\bm{b}_{n,n^{d}}\right)$.
From (\ref{eq: alltogether expand}), we observe that $m_{K_{n}}^{\psi}\in\mathcal{M}_{S}^{\psi}$,
with $K_{n}=\sum_{k=1}^{n^{d}}n_{k}$.

To obtain (\ref{eq: res 3}), we write 
\begin{align*}
\left\Vert \eta_{n}-m_{K_{n}}^{\psi}\right\Vert _{p} & =\left\Vert \sum_{k=1}^{n^{d}}\text{Gate}_{k}\left(\bm{x};\bm{\gamma}_{n}\right)f\left(\bm{y}|\bm{x}_{k}^{n}\right)-\sum_{k=1}^{n^{d}}\text{Gate}_{k}\left(\bm{x};\bm{\gamma}\right)h_{n_{k}}^{k}\left(\bm{y}|\bm{x}_{k}^{n}\right)\right\Vert _{p}\\
 & \le\sum_{k=1}^{n^{d}}\left\Vert \text{Gate}_{k}\left(\bm{x};\bm{\gamma}_{n}\right)\left[f\left(\bm{y}|\bm{x}_{k}^{n}\right)-h_{n_{k}}^{k}\left(\bm{y}|\bm{x}_{k}^{n}\right)\right]\right\Vert _{p}\text{.}
\end{align*}
By separability and Fubini's theorem, we then have 
\begin{align*}
\left\Vert \eta_{n}-m_{K_{n}}^{\psi}\right\Vert  & \le\sum_{k=1}^{n^{d}}\left\Vert \text{Gate}_{k}\left(\bm{x};\bm{\gamma}_{n}\right)\right\Vert _{p,\mathbb{X}}\left\Vert f\left(\bm{y}|\bm{x}_{k}^{n}\right)-h_{n_{k}}^{k}\left(\bm{y}|\bm{x}_{k}^{n}\right)\right\Vert _{p,\mathbb{Y}}\\
 & \le\sup_{k\in\left[n^{d}\right]}\left\Vert \text{Gate}_{k}\left(\bm{x};\bm{\gamma}_{n}\right)\right\Vert _{p,\mathbb{X}}\sum_{k=1}^{n^{d}}\left\Vert f\left(\bm{y}|\bm{x}_{k}^{n}\right)-h_{n_{k}}^{k}\left(\bm{y}|\bm{x}_{k}^{n}\right)\right\Vert _{p,\mathbb{Y}}\text{.}
\end{align*}

Let $C_{2}\left(n,\bm{\gamma}_{n}\right)=\sup_{k\in\left[n^{d}\right]}\left\Vert \text{Gate}_{k}\left(\bm{x};\bm{\gamma}_{n}\right)\right\Vert _{p,\mathbb{X}}$.
Then, we apply Lemma \ref{lem tin lem} $n^{d}$ times to establish
the existence of a constant $N_{2}\left(\epsilon,n,\bm{\gamma}_{n}\right)\in\mathbb{N}$,
such that for all $k\in\left[n^{d}\right]$ and $n_{k}\ge N_{2}\left(\epsilon,n,\bm{\gamma}_{n}\right)$,
\[
\left\Vert f\left(\bm{y}|\bm{x}_{k}^{n}\right)-h_{n_{k}}^{k}\left(\bm{y}|\bm{x}_{k}^{n}\right)\right\Vert _{p,\mathbb{Y}}\le\frac{\epsilon}{3C_{2}\left(n,\bm{\gamma}_{n}\right)n^{d}}\text{.}
\]
Thus, we have 
\[
\left\Vert \eta_{n}-m_{K_{n}}^{\psi}\right\Vert \le C_{2}\left(n,\bm{\gamma}_{n}\right)\times n^{d}\times\frac{\epsilon}{3C_{2}\left(n,\bm{\gamma}_{n}\right)n^{d}}=\frac{\epsilon}{3}\text{,}
\]
which completes our proof.

\subsection{Proof of Theorem \ref{thm main d=00003D00003D1}}

The proof is procedurally similar to that of Theorem \ref{thm main}
and thus we only seek to highlight the important differences. Firstly,
for any $\epsilon>0$, we approximate $f\left(\bm{y}|\bm{x}\right)$
by $\upsilon_{n}\left(\bm{x}|\bm{y}\right)$ of form (\ref{eq: discrete fine approx}),
with $d=1$. Result (\ref{eq: res 1}) implies uniform convergence,
in the sense that there exists an $N_{1}\left(\epsilon\right)\in\mathbb{N}$,
such that for all $n\ge N_{1}\left(\epsilon\right)$, 
\begin{equation}
\left\Vert f-\upsilon_{n}\right\Vert _{\mathcal{B}}<\frac{\epsilon}{3}\text{.}\label{eq: res 1 II}
\end{equation}

We now seek to approximate $\upsilon_{n}$ by $\eta_{n}$ of form
(\ref{eq: eta fun}), with $\bm{\gamma}_{n}=\bm{\gamma}_{l}$ for
some $l\in\mathbb{N}$. Upon application of Lemma \ref{lem almost uniform jiang},
it follows that for each $k\in\left[n^{d}\right]$ and $\varepsilon>0$,
there exists a measurable set $\mathbb{B}_{k}\left(\varepsilon\right)\subseteq\mathbb{X}$,
such that

\[
\lambda\left(\mathbb{B}_{k}\left(\varepsilon\right)\right)<\frac{\varepsilon}{n^{d}\lambda\left(\mathbb{Y}\right)}
\]
and 
\[
\left\Vert \mathrm{Gate}_{k}\left(\cdot;\bm{\gamma}_{l}\right)-\mathbf{1}_{\left\{ \bm{x}\in\mathbb{X}_{k}^{n}\right\} }\right\Vert _{\mathcal{B}\left(\mathbb{B}_{k}^{\complement}\left(\varepsilon\right)\right)}<\frac{\epsilon}{3}
\]
for all $l\ge M_{k}\left(\epsilon,n\right)$, for some $M_{k}\left(\epsilon,n\right)\in\mathbb{N}$.
Here, $\left(\cdot\right)^{\complement}$ is the set complement operator.

Since $f\in\mathcal{B}$, we have the bound $C\left(n\right)=\sum_{k=1}^{n^{d}}\left\Vert f\left(\bm{y}|\bm{x}_{k}^{n}\right)\right\Vert _{\mathcal{B}\left(\mathbb{Y}\right)}<\infty$.
Write $\mathbb{B}\left(\varepsilon\right)=\bigcup_{k=1}^{n^{d}}\mathbb{B}_{k}\left(\varepsilon\right)$.
Then, $\mathbb{B}^{\complement}\left(\varepsilon\right)=\bigcap_{k=1}^{n^{d}}\mathbb{B}_{k}^{\complement}\left(\varepsilon\right)$,
\[
\lambda\left(\mathbb{B}\left(\varepsilon\right)\right)\le\sum_{k=1}^{n^{d}}\lambda\left(\mathbb{B}_{k}^{\complement}\left(\varepsilon\right)\right)<\frac{\varepsilon}{\lambda\left(\mathbb{Y}\right)}\text{,}
\]
and 
\begin{align*}
\left\Vert \mathrm{Gate}_{k}\left(\cdot;\bm{\gamma}_{l}\right)-\mathbf{1}_{\left\{ \bm{x}\in\mathbb{X}_{k}^{n}\right\} }\right\Vert _{\mathcal{B}\left(\mathbb{B}^{\complement}\left(\varepsilon\right)\right)} & \le\min_{k\in\left[n^{d}\right]}\left\Vert \mathrm{Gate}_{k}\left(\cdot;\bm{\gamma}_{l}\right)-\mathbf{1}_{\left\{ \bm{x}\in\mathbb{X}_{k}^{n}\right\} }\right\Vert _{\mathcal{B}\left(\mathbb{B}_{k}^{\complement}\left(\varepsilon\right)\right)}<\frac{\epsilon}{3C\left(n\right)}\text{,}
\end{align*}
for all $l\ge M\left(\epsilon,n\right)=\max_{k\in\left[n^{d}\right]}M_{k}\left(\epsilon,n\right)$.
Here we use the fact that the supremum over some intersect of sets
is less than or equal to the minimum of the supremum over each individual
set.

Upon defining $\mathbb{C}\left(\varepsilon\right)=\mathbb{B}\left(\varepsilon\right)\times\mathbb{Y}\subset\mathbb{Z}$,
we observe that 
\[
\lambda\left(\mathbb{C}\left(\varepsilon\right)\right)=\lambda\left(\mathbb{B}\left(\varepsilon\right)\right)\lambda\left(\mathbb{Y}\right)\le\frac{\varepsilon}{\lambda\left(\mathbb{Y}\right)}\times\lambda\left(\mathbb{Y}\right)=\varepsilon\text{,}
\]
and $\mathbb{C}\left(\varepsilon\right)\subset\mathbb{B}\left(\varepsilon\right)\times\mathbb{Y}$.
Note also that 
\[
\left(\mathbb{B}\left(\varepsilon\right)\times\mathbb{Y}\right)^{\complement}=\mathbb{Z}\backslash\left(\mathbb{B}\left(\varepsilon\right)\times\mathbb{Y}\right)=\mathbb{B}^{\complement}\left(\varepsilon\right)\times\mathbb{Y}
\]
and 
\[
\mathbb{C}^{\complement}\left(\varepsilon\right)=\left(\mathbb{B}^{\complement}\left(\varepsilon\right)\times\mathbb{Y}\right)\cup\left(\mathbb{B}\left(\varepsilon\right)\times\mathbb{Y}^{\complement}\right)\cup\left(\mathbb{B}^{\complement}\left(\varepsilon\right)\times\mathbb{Y}^{\complement}\right)\text{.}
\]

It follows that

\[
\left\Vert \upsilon_{n}-\eta_{n}\right\Vert _{\mathcal{B\left(\mathbb{C}^{\complement}\left(\varepsilon\right)\right)}}\le\max\left\{ \left\Vert \upsilon_{n}-\eta_{n}\right\Vert _{\mathcal{B\left(\mathbb{B}^{\complement}\left(\varepsilon\right)\times\mathbb{Y}\right)}},\left\Vert \upsilon_{n}-\eta_{n}\right\Vert _{\mathcal{B\left(\mathbb{B}\left(\varepsilon\right)\times\mathbb{Y}^{\complement}\right)}},\left\Vert \upsilon_{n}-\eta_{n}\right\Vert _{\mathcal{B\left(\mathbb{B}^{\complement}\left(\varepsilon\right)\times\mathbb{Y}^{\complement}\right)}}\right\} \text{.}
\]
Since $\mathbb{B}\left(\varepsilon\right)\times\mathbb{Y}^{\complement}$
and $\mathbb{B}^{\complement}\left(\varepsilon\right)\times\mathbb{Y}^{\complement}$
are empty, via separability, we have 
\begin{align*}
\left\Vert \upsilon_{n}-\eta_{n}\right\Vert _{\mathcal{B\left(\mathbb{C}^{\complement}\left(\varepsilon\right)\right)}} & =\left\Vert \upsilon_{n}-\eta_{n}\right\Vert _{\mathcal{B\left(\mathbb{B}^{\complement}\left(\varepsilon\right)\times\mathbb{Y}\right)}}\\
 & =\sup_{\bm{z}\in\mathbb{B}^{\complement}\left(\varepsilon\right)\times\mathbb{Y}}\left|\sum_{k=1}^{n^{d}}\left[\mathbf{1}_{\left\{ \bm{x}\in\mathbb{X}_{k}^{n}\right\} }-\mathrm{Gate}_{k}\left(\bm{x};\bm{\gamma}_{l}\right)\right]f\left(\bm{y}|\bm{x}_{k}^{n}\right)\right|\\
 & \le\sup_{\bm{z}\in\mathbb{B}^{\complement}\left(\varepsilon\right)\times\mathbb{Y}}\sum_{k=1}^{n^{d}}\left|\mathbf{1}_{\left\{ \bm{x}\in\mathbb{X}_{k}^{n}\right\} }-\mathrm{Gate}_{k}\left(\bm{x};\bm{\gamma}_{l}\right)\right|f\left(\bm{y}|\bm{x}_{k}^{n}\right)\\
 & \le\sum_{k=1}^{n^{d}}\sup_{\bm{z}\in\mathbb{B}^{\complement}\left(\varepsilon\right)\times\mathbb{Y}}\left|\mathbf{1}_{\left\{ \bm{x}\in\mathbb{X}_{k}^{n}\right\} }-\mathrm{Gate}_{k}\left(\bm{x};\bm{\gamma}_{l}\right)\right|f\left(\bm{y}|\bm{x}_{k}^{n}\right)\\
 & =\sum_{k=1}^{n^{d}}\left\Vert \mathbf{1}_{\left\{ \bm{x}\in\mathbb{X}_{k}^{n}\right\} }-\mathrm{Gate}_{k}\left(\bm{x};\bm{\gamma}_{l}\right)\right\Vert _{\mathcal{B\left(\mathbb{B}^{\complement}\left(\varepsilon\right)\right)}}\left\Vert f\left(\bm{y}|\bm{x}_{k}^{n}\right)\right\Vert _{\mathcal{B\left(\mathbb{Y}\right)}}\\
 & \le\sup_{k\in\left[n\right]}\left\Vert \mathbf{1}_{\left\{ \bm{x}\in\mathbb{X}_{k}^{n}\right\} }-\mathrm{Gate}_{k}\left(\bm{x};\bm{\gamma}_{l}\right)\right\Vert _{\mathcal{B\left(\mathbb{B}^{\complement}\left(\varepsilon\right)\right)}}\sum_{k=1}^{n^{d}}\left\Vert f\left(\bm{y}|\bm{x}_{k}^{n}\right)\right\Vert _{\mathcal{B\left(\mathbb{Y}\right)}}\text{.}
\end{align*}
Recall that the $\sum_{k=1}^{n^{d}}\left\Vert f\left(\bm{y}|\bm{x}_{k}^{n}\right)\right\Vert _{\mathcal{B\left(\mathbb{Y}\right)}}=C\left(n\right)<\infty$
and that we can choose $l\ge M\left(\epsilon,n\right)$ so that 
\[
\sup_{k\in\left[n\right]}\left\Vert \mathbf{1}_{\left\{ \bm{x}\in\mathbb{X}_{k}^{n}\right\} }-\mathrm{Gate}_{k}\left(\bm{x};\bm{\gamma}_{l}\right)\right\Vert _{\mathcal{B\left(\mathbb{B}^{\complement}\left(\varepsilon\right)\right)}}<\frac{\epsilon}{3C\left(n\right)}\text{,}
\]
and thus 
\begin{equation}
\left\Vert \upsilon_{n}-\eta_{n}\right\Vert _{\mathcal{B\left(\mathbb{C}^{\complement}\left(\varepsilon\right)\right)}}<\frac{\epsilon}{3C\left(n\right)}\times C\left(n\right)=\frac{\epsilon}{3}\text{,}\label{eq: res 3 II}
\end{equation}
as required.

Finally, by noting that for each $k\in\left[n^{d}\right]$, both (\ref{eq: individual approx})
and $f\left(\cdot|\bm{x}_{k}^{n}\right)$ are continuous over $\mathbb{Y}$,
we apply Lemma \ref{lem tin lem} to obtain an $N_{2}\left(\epsilon,n,l\right)\in\mathbb{N}$,
such that for any $\epsilon>0$ and $n_{k}\ge N_{2}\left(\epsilon,n,l\right)$,
we have
\[
\left\Vert f\left(\cdot|\bm{x}_{k}^{n}\right)-h_{n_{k}}^{k}\left(\cdot|\bm{x}_{k}^{n}\right)\right\Vert _{\mathcal{B}\left(\mathbb{Y}\right)}<\frac{\epsilon}{3M_{1}n}\text{.}
\]
Here $M_{1}=\sup_{k\in\left[n^{d}\right]}\left\Vert \mathrm{Gate}_{k}\left(\cdot;\bm{\gamma}_{l}\right)\right\Vert _{\mathcal{B\left(\mathbb{X}\right)}}<\infty$,
since $\mathrm{Gate}_{k}\left(\bm{x};\bm{\gamma}_{l}\right)$ is continuous
in $\bm{x}$, and $\mathbb{X}$ is compact. Therefore, for all $K_{n}=\sum_{k=1}^{n^{d}}n_{k}$,
$n_{k}\ge N_{2}\left(\epsilon,n,l\right)$,
\begin{align}
\left\Vert \eta_{n}-m_{K_{n}}^{\psi}\right\Vert _{\mathcal{B}} & \le\sup_{k\in\left[n^{d}\right]}\left\Vert \mathrm{Gate}_{k}\left(\bm{x};\bm{\gamma}_{l}\right)\right\Vert _{\mathcal{B\left(\mathbb{X}\right)}}\sum_{k=1}^{n^{d}}\left\Vert f\left(\cdot|\bm{x}_{k}^{n}\right)-h_{n_{k}}^{k}\left(\cdot|\bm{x}_{k}^{n}\right)\right\Vert _{\mathcal{B}\left(\mathbb{Y}\right)}\nonumber \\
 & =M_{1}\times n^{d}\times\frac{\epsilon}{3M_{1}n^{d}}=\frac{\epsilon}{3}.\label{eq: res 3 II-1}
\end{align}

In summary, via \eqref{eq: res 1 II}, \eqref{eq: res 3 II}, and
\eqref{eq: res 3 II-1}, for each $\epsilon>0$, for any $\varepsilon>0$,
there exists a $\mathbb{C}\left(\varepsilon\right)\subset\mathbb{Z}$
and constants $N_{1}\left(\epsilon\right),M\left(\epsilon,n\right),N_{2}\left(\epsilon,n,l\right)\in\mathbb{N}$,
such that for all $K_{n}=\sum_{k=1}^{n^{d}}n_{k}$, with $n_{k}\ge N_{2}\left(\epsilon,n,l\right)$,
$l\ge M\left(\epsilon,n\right)$, and $n\ge N_{1}\left(\epsilon\right)$,
it follows that $\lambda\left(\mathbb{C}\left(\varepsilon\right)\right)<\varepsilon$,
and 
\begin{align*}
\left\Vert f-m_{K_{n}}^{\psi}\right\Vert _{\mathcal{B\left(\mathbb{C}^{\complement}\left(\varepsilon\right)\right)}} & \le\left\Vert f-\upsilon_{n}\right\Vert _{\mathcal{B\left(\mathbb{C}^{\complement}\left(\varepsilon\right)\right)}}+\left\Vert \upsilon_{n}-\eta_{n}\right\Vert _{\mathcal{B\left(\mathbb{C}^{\complement}\left(\varepsilon\right)\right)}}+\left\Vert \eta_{n}-m_{K_{n}}^{\psi}\right\Vert _{\mathcal{B\left(\mathbb{C}^{\complement}\left(\varepsilon\right)\right)}}\\
 & \le\left\Vert f-\upsilon_{n}\right\Vert _{\mathcal{B}}+\left\Vert \upsilon_{n}-\eta_{n}\right\Vert _{\mathcal{B\left(\mathbb{C}^{\complement}\left(\varepsilon\right)\right)}}+\left\Vert \eta_{n}-m_{K_{n}}^{\psi}\right\Vert _{\mathcal{B}}\\
 & <3\times\frac{\epsilon}{3}=\epsilon\text{.}
\end{align*}
This completes the proof.

\section{\label{sec:Proofs-of-lemmas}Proofs of lemmas}

\subsection{Proof of Lemma \ref{lem equiv}}

We firstly prove that any gating vector from $\mathcal{G}_{S}^{K}$
can be equivalently represented as an element of $\mathcal{G}_{G}^{K}$.
For any $\bm{x}\in\mathbb{R}^{d}$, $d\in\mathbb{N}$, $k\in\left[K\right]$,
$a_{k}\in\mathbb{R}$, $\bm{b}_{k}\in\mathbb{R}^{d}$, and $K\in\mathbb{N}$,
choose $\bm{\nu}_{k}=\bm{b}_{k}$, $\tau_{k}=a_{k}+\bm{b}_{k}^{\top}\bm{b}_{k}/2$
and 
\[
\pi_{k}=\exp\left(\tau_{k}\right)/\sum_{l=1}^{K}\exp\left(\tau_{l}\right)\text{.}
\]
This implies that $\sum_{l=1}^{K}\pi_{l}=1$, $\pi_{l}>0$, for all
$l\in\left[K\right]$, and 
\begin{align*}
\frac{\exp\left(a_{k}+\bm{b}_{k}^{\top}\bm{x}\right)}{\sum_{l=1}^{K}\exp\left(a_{k}+\bm{b}_{k}^{\top}\bm{x}\right)} & =\frac{\exp\left(\tau_{k}-\bm{v}_{k}^{\top}\bm{v}_{k}/2+\bm{v}_{k}^{\top}\bm{x}\right)}{\sum_{l=1}^{K}\exp\left(\tau_{l}-\bm{v}_{l}^{\top}\bm{v}_{l}/2+\bm{v}_{l}^{\top}\bm{x}\right)}\\
 & =\frac{\exp\left(\tau_{k}\right)\exp\left(-\left(\bm{x}-\bm{\nu}_{k}\right)^{\top}\left(\bm{x}-\bm{\nu}_{k}\right)/2\right)}{\sum_{l=1}^{K}\exp\left(\tau_{l}\right)\exp\left(-\left(\bm{x}-\bm{\nu}_{l}\right)^{\top}\left(\bm{x}-\bm{\nu}_{l}\right)/2\right)}\\
 & =\frac{\pi_{k}\left(2\pi\right)^{-d/2}\exp\left(-\left(\bm{x}-\bm{\nu}_{k}\right)^{\top}\left(\bm{x}-\bm{\nu}_{k}\right)/2\right)}{\sum_{l=1}^{K}\pi_{l}\left(2\pi\right)^{-d/2}\exp\left(-\left(\bm{x}-\bm{\nu}_{l}\right)^{\top}\left(\bm{x}-\bm{\nu}_{l}\right)/2\right)}\\
 & =\frac{\pi_{k}\phi\left(\bm{x};\bm{\nu}_{k},\mathbf{I}\right)}{\sum_{l=1}^{K}\pi_{l}\phi\left(\bm{x};\bm{\nu}_{l},\mathbf{I}\right)}\text{,}
\end{align*}
where $\mathbf{I}$ is the identity matrix of appropriate size. This
proves that $\mathcal{G}_{S}^{K}\subset\mathcal{G}_{G}^{K}$.

Next, to show that $\mathcal{G}_{E}^{K}\subset\mathcal{G}_{S}^{K}$,
we write 
\begin{align*}
 & \frac{\pi_{k}\phi\left(\bm{x};\bm{\nu}_{k},\bm{\Sigma}\right)}{\sum_{l=1}^{K}\pi_{l}\phi\left(\bm{x};\bm{\nu}_{l},\bm{\Sigma}\right)}\\
= & \frac{\pi_{k}\left|2\pi\bm{\Sigma}\right|^{-1/2}\exp\left(-\left(\bm{x}-\bm{\nu}_{k}\right)^{\top}\bm{\Sigma}^{-1}\left(\bm{x}-\bm{\nu}_{k}\right)/2\right)}{\sum_{l=1}^{K}\pi_{l}\left|2\pi\bm{\Sigma}\right|^{-1/2}\exp\left(-\left(\bm{x}-\bm{\nu}_{l}\right)^{\top}\bm{\Sigma}^{-1}\left(\bm{x}-\bm{\nu}_{l}\right)/2\right)}\\
= & \frac{1}{\sum_{l=1}^{K}\exp\left(-\log\left(\pi_{l}^{-2}/\pi_{k}^{-2}\right)/2-\left(\bm{x}-\bm{\nu}_{l}\right)^{\top}\bm{\Sigma}^{-1}\left(\bm{x}-\bm{\nu}_{l}\right)/2-\left(\bm{x}-\bm{\nu}_{k}\right)^{\top}\bm{\Sigma}^{-1}\left(\bm{x}-\bm{\nu}_{k}\right)/2\right)}\text{,}
\end{align*}
and note that 
\begin{align*}
 & \left(\bm{x}-\bm{\nu}_{l}\right)^{\top}\bm{\Sigma}^{-1}\left(\bm{x}-\bm{\nu}_{l}\right)-\left(\bm{x}-\bm{\nu}_{k}\right)^{\top}\bm{\Sigma}^{-1}\left(\bm{x}-\bm{\nu}_{k}\right)\\
= & -2\left(\bm{\nu}_{l}-\bm{\nu}_{k}\right)^{\top}\bm{\Sigma}^{-1}\bm{x}+\left(\bm{\nu}_{l}+\bm{\nu}_{k}\right)^{\top}\bm{\Sigma}^{-1}\left(\bm{\nu}_{l}-\bm{\nu}_{k}\right)\text{.}
\end{align*}
Thus, we have 
\begin{align*}
 & \frac{\pi_{k}\phi\left(\bm{x};\bm{\nu}_{k},\bm{\Sigma}\right)}{\sum_{l=1}^{K}\pi_{l}\phi\left(\bm{x};\bm{\nu}_{l},\bm{\Sigma}\right)}\\
= & \frac{1}{\sum_{l=1}^{K}\exp\left(-\log\left(\pi_{l}^{-2}/\pi_{k}^{-2}\right)/2-\left(\bm{\nu}_{l}+\bm{\nu}_{k}\right)^{\top}\bm{\Sigma}^{-1}\left(\bm{\nu}_{l}-\bm{\nu}_{k}\right)/2-\left(\bm{\nu}_{l}-\bm{\nu}_{k}\right)^{\top}\bm{\Sigma}^{-1}\bm{x}\right)}\text{.}
\end{align*}

Next, notice that we can write 
\[
\frac{\exp\left(a_{k}+\bm{b}_{k}^{\top}\bm{x}\right)}{\sum_{l=1}^{K}\exp\left(a_{l}+\bm{b}_{l}^{\top}\bm{x}\right)}=\frac{1}{\sum_{l=1}^{K}\exp\left(\alpha_{l}+\bm{\beta}_{l}^{\top}\bm{x}\right)}\text{,}
\]
where $\alpha_{l}=a_{l}-a_{k}$ and $\bm{\beta}_{l}=\bm{\beta}_{l}-\bm{\beta}_{k}$.
We now choose $a_{k}$ and $\bm{b}_{k}$, such that for every $l\in\left[K\right]$,
\begin{align*}
\alpha_{l} & =a_{l}-a_{k}=-\frac{1}{2}\log\left(\frac{\pi_{l}^{-2}}{\pi_{k}^{-2}}\right)-\frac{1}{2}\left(\bm{\nu}_{l}^{\top}\bm{\Sigma}^{-1}\bm{\nu}_{l}-\bm{\nu}_{k}^{\top}\bm{\Sigma}^{-1}\bm{\nu}_{k}\right)\text{,}
\end{align*}
and 
\[
\bm{\beta}_{l}=\bm{\beta}_{l}-\bm{\beta}_{k}=\bm{\nu}_{l}^{\top}\bm{\Sigma}^{-1}-\bm{\nu}_{k}^{\top}\bm{\Sigma}^{-1}\text{.}
\]
To complete the proof, we choose 
\[
a_{k}=\log\left(\pi_{k}\right)-\frac{1}{2}\bm{\nu}_{k}^{\top}\bm{\Sigma}^{-1}\bm{\nu}_{k}
\]
and $b_{k}=\bm{\nu}_{k}^{\top}\bm{\Sigma}^{-1}$, for each $k\in\left[K\right]$.

\subsection{Proof of Lemma \ref{lem almost uniform jiang}}

For $l\in\left[0,\infty\right)$, write 
\[
\text{Gate}_{k}\left(x,l\right)=\frac{\exp\left(\left[x-c_{k}\right]lk\right)}{\sum_{i=1}^{n}\exp\left(\left[x-c_{i}\right]li\right)}\text{,}
\]
where $x\in\mathbb{X}=\left[0,1\right]$, and $c_{k}=\left(k-1\right)/\left(2k\right)$.
We identify that $\mathbf{Gate}=\left(\text{Gate}_{k}\left(x,l\right)\right)_{k\in\left[n\right]}$
belongs to the class $\mathcal{G}_{S}^{n}$. The proof of the Section
4 Proposition from \citet{Jiang1999} reveals that for all $k\in\left[n\right]$,
\[
\text{Gate}_{k}\left(x,l\right)\rightarrow\mathbf{1}_{\left\{ x\in\mathbb{I}_{k}^{n}\right\} }
\]
almost everywhere in $\lambda$, as $l\rightarrow\infty$. The result
then follows via an application of Egorov's theorem (cf. \citealt[Thm. 2.33]{Folland:1999aa}).

\section{Summary and conclusions\label{sec:Summary-and-conclusions}}

Using recent results mixture model approximation results \citet{Nguyen:2020aa}
and \citet{Nguyen:2020ab}, and the indicator approximation theorem
of \citet{Jiang1999} (cf. Section \ref{sec:Technical-lemmas}), we
have proved two approximation theorems (Theorems \ref{thm main} and
\ref{thm main d=00003D00003D1}) regarding the class of soft-max gated
MoE models with experts arising from arbitrary location-scale families
of conditional density functions. Via an equivalence result (Lemma
\ref{lem equiv}), the results of Theorems \ref{thm main} and \ref{thm main d=00003D00003D1}
also extend to the setting of Gaussian gated MoE models (Corollary
\ref{cor: Gaussian gate}), which can be seen as a generalization
of the soft-max gated MoE models.

Although we explicitly make the assumption that $\mathbb{X}=\left[0,1\right]^{d}$,
for the sake of mathematical argument (so that we can make direct
use of Lemma \ref{Lem Jiang}), a simple shift-and-scale argument
can be used to generalize our result to cases where $\mathbb{X}$
is any generic compact domain. The compactness assumption regarding
the input domain is common in the MoE and mixture of regression models
literature, as per the works of \citet{Jiang1999a}, \citet{Jiang1999},
\citet{Norets2010}, \citet{Montuelle2014}, \citet{Pelenis2014},
\citet{Devijver2015An-ell1-oracle-}, and \citet{Devijver2015Finite-mixture-}.

The assumption permits the application of the result to the settings
where the inputs $\bm{X}$ is assumed to be non-random design vectors
that take value on some compact set $\mathbb{X}$. This is often the
case when there is only a finite number of possible design vector
elements for which $\bm{X}$ can take. Otherwise, the assumption also
permits the scenario where $\bm{X}$ is some random element with compactly
supported distribution, such as uniformly distributed, or beta distributed
inputs. Unfortunately, the case of random $\bm{X}$ over an unbounded
domain (e.g., if $\bm{X}$ has multivariate Gaussian distribution)
is not covered under our framework. An extension to such cases would
require a more general version of Lemma \ref{Lem Jiang}, which we
believe is a nontrivial direction for future work.

Like the input, we also assume that the output domain is restricted
to a compact set $\mathbb{Y}$. However, the output domain of the
approximating class of MoE models is unrestricted to $\mathbb{Y}$
and thus the functions (i.e., we allow $\psi$ to be a PDF over $\mathbb{R}^{q}$).
The restrictions placed on $\mathbb{Y}$ is also common in the mixture
approximation literature, as per the works of \citet{ZeeviMeir1997},
\citet{LiBarron1999}, and \citet{RakhlinPanchenkoMukherjee2005},
and is also often made in the context of nonparametric regression
(see, e.g., \citealp{Gyorfi2002} and \citealp{Cucker2007Learning-Theory}).
Here, our use of the compactness of $\mathbb{Y}$ is to bound the
integral of $v_{n}$, in \eqref{eq: control v_n}. A more nuanced
approach, such as via the use of a generalized Lebesgue spaces (see
e.g., \citealp{Castillo2010An-Introductory} and \citealp{Cruze-Uribe2013Variable-Lebesg}),
may lead to result for unbounded $\mathbb{Y}$. This is another exciting
future direction of our research program. 

A trivial modification to the proof of Lemma \ref{lem tin lem} allows
us to replace the assumption that $f$ is a PDF with a sub-PDF assumption
(i.e., $\int_{\mathbb{Y}}f\text{d}\lambda\le1$), instead. This in
turn permits us to replace the assumption that $f\left(\cdot|\bm{x}\right)$
is a conditional PDF in Theorems \ref{thm main} and \ref{thm main d=00003D00003D1}
with sub-PDF assumptions as well (i.e., for each $\bm{x}\in\mathbb{X}$,
$\int_{\mathbb{Y}}f\left(\bm{y}|\bm{x}\right)\text{d}\lambda\left(\bm{y}\right)\le1$).
Thus, in this modified form, we have a useful interpretation for situations
when the input $\bm{Y}$ is unbounded. That is, when $\bm{Y}$ is
unbounded, we can say that the conditional PDF $f$ can be arbitrarily
well approximated in $\mathcal{L}_{p}$ norm by a sequence $\left\{ m_{K}^{\psi}\right\} _{K\in\mathbb{N}}$
of either soft-max or Gaussian gated MoEs over any compact subdomain
$\mathbb{Y}$ of the unbounded domain of $\bm{Y}$. Thus, although
we cannot provide guarantees of the entire domain of $\bm{Y}$, we
are able to guarantee arbitrary approximate fidelity over any arbitrarily
large compact subdomain. This is a useful result in practice since
one is often not interested in the entire domain of $\bm{Y}$, but
only on some subdomain where the probability of $\bm{Y}$ is concentrated.
This version of the result resembles traditional denseness results
in approximation theory, such as those of \citet[Ch. 20]{CheneyLight2000}.

Finally, our results can be directly applied to provide approximation
guarantees for a large number of currently used models in applied
statistics and machine learning research. Particularly, our approximation
guarantees are applicable to the recent MoE models of \citet{Ingrassia2012},
\citet{Chamroukhi2013a}, \citet{Ingrassia2014}, \citet{Chamroukhi2016},
\citet{Nguyen2016}, \citet{Deleforge2015}, \citet{Deleforge2015a},
\citet{Kalliovirta2016}, and \citet{Perthame2018}, among many others.
Here, we may guarantee that the underlying data generating processes,
if satisfying our assumptions, can be adequately well approximated
by sufficiently complex forms of the models considered in each of
the aforementioned work. 

The rate and manner of which good approximation can be achieved as
a function of the number of experts $K$ and the sample size is a
currently active research area, with pioneering work conducted in
\citet{Cohen2012Conditional-den} and \citet{Montuelle2014}. More
recent results in this direction appear in \citet{Nguyen2020An-l1-oracle-in}, \citet{Nguyen2021A-non-asymptoti}, and \citet{nguyen2021nonBLoMPE}.

\section*{List of abbreviations}

\begin{description}
\item [{MoE}] Mixture of experts
\item [{PDF}] Probability density function
\end{description}

\section*{Acknowledgements}

Hien Duy Nguyen and Geoffrey John McLachlan are funded by Australian Research
Council grants: DP180101192 and IC170100035. TrungTin Nguyen is supported by
a \textquotedblleft Contrat doctoral\textquotedblright{} from the
French Ministry of Higher Educationand Research. Faicel Chamroukhi
is granted by the French National Research Agency (ANR) grant \href{https://anr.fr/en/funded-projects-and-impact/funded-projects/project/funded/project/b2d9d3668f92a3b9fbbf7866072501ef-f004f5ad27/?tx_anrprojects\_funded\%5Bcontroller\%5D=Funded&cHash=895d4c3a16ad6a0902e6515eb65fba37}{SMILES ANR-18-CE40-0014}. The authors also thank the Editor and Reviewer,
whose careful and considerate comments lead to improvements in of
text.

\section*{Declarations}

\textbf{Funding:} HDN and GJM are funded by Australian Research Council
grants: DP180101192 and IC170100035. FC is funded by ANR grant: \href{https://anr.fr/en/funded-projects-and-impact/funded-projects/project/funded/project/b2d9d3668f92a3b9fbbf7866072501ef-f004f5ad27/?tx_anrprojects\_funded\%5Bcontroller\%5D=Funded&cHash=895d4c3a16ad6a0902e6515eb65fba37}{SMILES ANR-18-CE40-0014}.\\
 \textbf{Conflicts of interest:} None.\\
 \textbf{Availability of data and material: }Not applicable.\\
  \textbf{Authors' contributions: }All authors contributed equally to the exposition and to the mathematical derivations.\\
 \textbf{Code availability:} Not applicable.

 \bibliographystyle{apalike2}
\bibliography{20201118_MASTERBIB}

\end{document}